\documentclass[12pt]{article}
\usepackage{amsmath}
\newtheorem{theo}{Theorem}[section]
\def\qed{\hfill \rule{4pt}{7pt}}
\def\pf{\noindent {\it{Proof.} \hskip 2pt}}
\newcommand{\A}{\mathcal{A}}
\newcommand{\US}{\mathcal{US}}
\newcommand{\NCSym}{\mathrm{NCSym}}

\begin{document}
\begin{center}
{\Large\bf  A Bijection between
 Atomic Partitions and

  Unsplitable
Partitions}
\end{center}

\begin{center}
William Y.C. Chen$^1$, Teresa X.S. Li$^2$  and David G.L. Wang$^3$

$^{1,2}$Center for Combinatorics, LPMC-TJKLC\\
Nankai University, Tianjin 300071, P.R. China

$^3$Beijing International Center for Mathematical Research\\
Peking University, Beijing 100871, P.R. China

$^1$chen@nankai.edu.cn, $^2$lxs@cfc.nankai.edu.cn,
$^3$wgl@math.pku.edu.cn
\end{center}

\begin{abstract}
In the study of the algebra
 $\NCSym$ of symmetric functions in
 noncommutative variables,
 Bergeron and Zabrocki found a free generating
  set consisting of
   power sum symmetric functions indexed by
 atomic partitions. On the other hand, Bergeron, Reutenauer,
Rosas, and Zabrocki studied
 another free generating set of
$\NCSym$ consisting of monomial symmetric functions indexed by
unsplitable partitions. Can and Sagan raised the question of
finding a bijection between atomic partitions
 and unsplitable partitions.
In this paper, we provide such a bijection.
\end{abstract}

\section{Introduction}

In their study of the algebra $\NCSym$ of symmetric
functions in noncommutative variables, Rosas and Sagan
\cite{Rosas-Sagan} introduced a vector space with a  basis
\[
\{p_{\pi}\,|\,\pi\textrm{ is a set partition}\},
\]
where $p_\pi$ is the power sum symmetric function in
noncommutative variables.
 Bergeron, Hohlweg, Rosas, and Zabrocki
\cite{BHRZ06} obtained the following formula
\[
p_{\pi|\sigma}=p_{\pi}\,p_{\sigma},
\]
where $\pi|\sigma$ denotes the slash product of $\pi$ and
$\sigma$. It follows that, as an algebra, $\NCSym$ is freely
generated by $p_{\pi}$ with $\pi$ atomic, see Bergeron and
Zabrocki \cite{BZ05X}. It should be noted that  Wolf \cite{Wolf}
showed that  $\NCSym$  is freely generated by another basis. A
combinatorial characterization of the generating set of Wolf has
been found by Bergeron, Reutenauer, Rosas, and Zabrocki
\cite{BERGERON-R-R-Z}. More precisely, they introduced the notion
of unsplitable partitions and proved
 that the generating set of
Wolf can be described as the set of monomial
symmetric functions in noncommutative
variables indexed by unsplitable partitions.

Let $[n]$ denote the set $\{1,2,\ldots,n\}$. Taking the degree into
account, one sees that the number of atomic partitions of $[n]$
equals the number of unsplitable partitions of $[n]$. Recently,
Can and Sagan~\cite{Can-Sagan} raised the question of finding a
combinatorial proof of this fact. The objective of this paper is
to present such a proof.

%\begin{example}
%For $n=4$, there are $6$ atomic partitions:
%\begin{xalignat*}{3}
%&\bigl\{\{1,2,3,4\}\bigr\},
%&& \bigl\{\{1,2,4\},\,\{3\}\bigr\},
%&& \bigl\{\{1,3,4\},\,\{2\}\bigr\},\\[4pt]
%&\bigl\{\{1,3\},\,\{2,4\}\bigr\},
%&&\bigl\{\{1,4\},\,\{2,3\}\bigr\},
%&&\bigl\{\{1,4\},\,\{2\},\,\{3\}\bigr\}.
%\end{xalignat*}
%And there are also $6$ unsplitable partitions, namely,
%\begin{xalignat*}{3}
%&\bigl\{\{1,3\},\,\{2\},\,\{4\}\bigr\},
%&&\bigl\{\{1\},\,\{2,3,4\}\bigr\},
%&&\bigl\{\{1\},\,\{2,3\},\,\{4\}\bigr\},\\[4pt]
%&\bigl\{\{1\},\,\{2,4\},\,\{3\}\bigr\},
%&&\bigl\{\{1\},\,\{2\},\,\{3,4\}\bigr\},
%&&\bigl\{\{1\},\,\{2\},\,\{3\},\,\{4\}\bigr\}.
%\end{xalignat*}
%\end{example}

\section{The bijection}

In this section we construct a bijection between the set of atomic
partitions of $[n]$ and the set of unsplitable partitions of
$[n]$.

Let us begin with an overview of terminology. Let $X$ be an
totally ordered set. A partition $\pi$ of $X$ is a family
$\{B_1,B_2,\ldots,B_k\}$ of disjoint nonempty subsets of $X$
 whose
union is  $X$. The subsets $B_i$
 are called blocks of $\pi$.
Without loss of generality, we may assume that the blocks of a
partition are arranged in the increasing order of
 their minimal
elements, and that the elements in each block are written in
increasing order.

 Let $\pi$ be a partition of $X$ and $S\subseteq
X$. We say that $\sigma$ is the restriction of $\pi$ on $S$,
denoted by $\sigma=\pi_{S}$, if $\sigma$ is a partition of
 $S$ such that any
two elements lie in the same block of $\sigma$
 if and only if they
are in the same block of $\pi$. In other words, $\pi_S$ is
obtained from $\pi$ by removing all elements
 that do not belong to
$S$. For two positive integers $i$ and $j$ with $i<j$, we use
 $[i,j]$ to denote the set $\{i,i+1,\ldots,j\}$.
  For example, if
\begin{equation}\label{eg_pi}
\pi=\bigl\{\{1,3,5,6\},\,\{2,7,9\},\,\{4,8,10\}\bigr\},
\end{equation}
then
\begin{equation}\label{eg_restriction}
\pi_{[5,10]}=\bigl\{\{5,6\},\,\{7,9\},\,\{8,10\}\bigr\}.
\end{equation}

Let $\Pi_n$ be the set of partitions of $[n]$. Assume that
\begin{equation*}
\pi=\{B_1,B_2,\ldots,B_k\}\in\Pi_m, \qquad \sigma
=\{C_1,C_2,\ldots,C_l\}\in\Pi_n.
\end{equation*}
The slash product of $\pi$ and $\sigma$, denoted by
$\pi|\sigma$, is defined to be the partition obtained  by
joining the
 blocks of $\pi$ and the blocks of the
partition
\[
\sigma+m=\{C_1+m,\,C_2+m,\,\ldots,\,C_l+m\},
\]
that is,
\[
\pi|\sigma=
\{B_1,\,B_2,\,\ldots,\,B_k,\,C_1+m,\,C_2+m,\,\ldots,\,C_l+m\},
\]
where $C_i+m$ denotes the block obtained by adding $m$ to each
element in $C_i$. It can be seen that
$\pi|\sigma\in\Pi_{m+n}$. A partition $\pi$ is said to be
atomic if there are no nonempty partitions $\sigma$ and $\tau$
such that $\pi=\sigma|\tau$. Let $\A_n$ be the set of atomic
partitions of $[n]$. For example, for $n=3$ there are two atomic
partitions $\bigl\{\{1,3\},\,\{2\}\bigr\}$ and
$\bigl\{\{1,2,3\}\bigr\}$.

The split product of $\pi$ and $\sigma$, denoted by
$\pi\circ\sigma$, is given by
\[
\pi\circ\sigma=\begin{cases} \{B_1\cup (C_1+m),\,\ldots,\,B_k\cup
(C_k+m),\, C_{k+1}+m,\,\ldots,\,C_{l}+m\},
&\textrm{if } k\leq l; \\[5pt]
\{B_1\cup (C_1+m),\,\ldots,\,B_l\cup (C_l+m),\,B_{l+1},\,
\ldots,\,B_{k}\}, &\textrm{if } k>l.
\end{cases}
\]
Clearly, $\pi\circ\sigma\in\Pi_{m+n}$.
 A partition is said to be
splitable  if it is  the split product of two nonempty partitions.
Otherwise, it is said to be unsplitable.
 Denote by $\US_n$ the set
of unsplitable partitions of $[n]$. For example, for $n=3$ there
are two unsplitable partitions $\bigl\{\{1\},\,\{2,3\}\bigr\}$ and
$\bigl\{\{1\},\,\{2\},\,\{3\}\bigr\}$.

To describe our bijection, we first notice that it is possible for
a partition to be atomic and unsplitable at the same time. For example,
the partition
\[
\bigl\{\{1,3,7\},\,\{2,6\},\,\{4,5,8\}\bigr\}
\]
is both atomic and unsplitable. Our bijection will be concerned
with atomic partitions that are splitable and unsplitable
partitions that are not atomic.
 In other words, we shall
establish a bijection
\[
\varphi\colon \A_n\backslash\US_n\longrightarrow
\US_n\backslash\A_n.
\]

For the sake of presentation, let us introduce a notation. Let
$X=\{x_1,x_2,\ldots,x_n\}$ under the assumption that
$x_1<x_2<\cdots<x_n$.
Assume that $\pi=\{B_1,B_2,\ldots,B_k\}$ is
a partition of $X$. Let $r$ be the largest integer
$j$ such that
\[
B_j\cup B_{j+1}\cup\cdots\cup B_k=\{x_t,\,x_{t+1},\,\ldots,\,x_n\}
\]
for some $t$. The existence of such an integer $r$ is evident.
We define
\[
R(\pi)=\{B_r,\,B_{r+1},\,\ldots,\,B_k\}.
\]
Given the partition
$\pi=\bigl\{\{1,3,5,6\},\,\{2,7,9\},\,\{4,8,10\}\bigr\}$
as in \eqref{eg_pi}, we have
\begin{equation}\label{eg_raf}
R(\pi_{[5,10]})=\bigl\{\{7,9\},\{8,10\}\bigr\}.
\end{equation}
In the above notation, we see that  $\pi$ is atomic if and only if
$\pi=R(\pi)$.

We are now ready to present the map $\varphi$. Suppose that
$\pi=\{B_1,B_2,\ldots,B_k\}\in\A_n\backslash\US_n$. It
consists of three steps.
\begin{itemize}
\item[Step 1.] Let $i$ be the smallest element in~$B_1$ such that
$\pi=\pi_{[i-1]}\circ (\pi_{[i,\,n]}-i+1)$. The existence of the
element $i$  is guaranteed by the condition that $\pi$ is
splitable.
  \item[Step 2.] Let $j$ be the
smallest element in the underlying set of the
partition~$R(\pi_{[i,\,n]})$.
We see that $2\le i\le j\le n$ and
$R(\pi_{[i,\,n]})=\pi_{[j,\,n]}$.
 \item[Step 3.] Set
$\varphi(\pi)$ to be the partition $\pi_{[j-1]}\bigm|(\pi_{[j,\,n]}-j+1)$.
\end{itemize}

\begin{theo}
The map $\varphi$ is
 a bijection from $\A_n\backslash\US_n$
to $\US_n\backslash\A_n$.
\end{theo}

\pf First, we claim
 that $\varphi(\pi)\in\US_n\backslash\A_n$.
Since $2\le j\le n$, both $\pi_{[j-1]}$ and $\pi_{[j,\,n]}$ are
nonempty partitions. This implies  that $\varphi(\pi)\not\in\A_n$.

We next proceed to show that $\varphi(\pi)$ is unsplitable. To
this end, let
\begin{equation*}
\pi_{[j-1]} =\{C_1,C_2,\ldots,C_s\}, \qquad \pi_{[j,n]}
=\{D_1,D_2,\ldots,D_t\}.
\end{equation*}
Then
\[
\varphi(\pi)=
\{C_1,\,C_2,\,\ldots,\,C_s,\,D_1,\,D_2,\,\ldots,\,D_t\}.
\]
Suppose to the contrary that $\varphi(\pi)$
 is splitable, namely,
there exists an element $l\in C_1$ such that
\[
\varphi(\pi)=\varphi(\pi)_{[l-1]}\circ(\varphi(\pi)_{[l,\,n]}-l+1).
\]
Since $n$ belongs to some block $D_h$, by the definition of the
split product, we deduce that
\begin{equation}\label{1}
C_p\cap[l,\,n]\neq\emptyset,\quad\textrm{for each } 1\leq p\leq s.
\end{equation}
By the choice of $i$, we find that $l\geq i$. Recall that
$\pi=\{B_1,B_2,\ldots,B_k\}$. By the definition of $\pi_{[j,n]}$,
we may assume that the block $D_1$ of $\pi_{[j,n]}$ is contained
in some block $B_r$ of $\pi$. If $D_1=B_r$, then the smallest
element of $B_r$ is $j$. Therefore all elements in
$B_{r+1},B_{r+2},\ldots,B_k$ are larger than $j$. Now, by the
choice  of $j$, we deduce that
\[
B_r\cup B_{r+1}\cup\cdots\cup B_k=[j,n].
\]
Consequently,
\[
\pi=\pi_{[j-1]}\bigm|(\pi_{[j,\,n]}-j+1),
\]
which contradicts the assumption that $\pi$ is atomic.
Hence we have
$D_1\ne B_r$, and  so  $C_r=B_r\backslash
D_1\ne\emptyset$. Since $D_1$ is a block of the partition
$\pi_{[i,n]}$, it consists of
 all the elements in $B_r$ that are larger than or equal to
$i$. In other words,
each element in $C_r$ is less than $i$.
This yields that $C_r\cap[l,\,n]=\emptyset$,
a contradiction to \eqref{1}. Thus we have proved the claim
that  $\varphi(\pi)\in\US_n\backslash\A_n$.

We now define a map
\[
\psi\colon\US_n\backslash\A_n\longrightarrow
\A_n\backslash\US_n,
\]
and
shall show that $\psi$ is the inverse of $\varphi$.
Let $\sigma=\{B_1,B_2,\ldots,B_k\}\in\US_n\backslash\A_n$.
\begin{itemize}
\item[Step 1.]
Let $j$ be the smallest element in the
underlying set of the partition $R(\sigma)$.
\item[Step 2.]
Let $B_r$ be the first block in the partition $R(\sigma)$.
We consider two cases.
\begin{itemize}
\item[Case 1.]
If $\sigma_{[j-1]}$ is unsplitable, then set
\[
\psi(\sigma)=\sigma_{[j-1]}\circ(\sigma_{[j,\,n]}-j+1).
\]
\item[Case 2.]
If $\sigma_{[j-1]}$ is splitable, then choose $i$  to be the smallest
element in $B_1$ such that
\begin{equation}\label{j-1}
\sigma_{[j-1]}=\sigma_{[i-1]}\circ(\sigma_{[i,\,j-1]}-i+1).
\end{equation}
Let $q=\min\{l\,|\,B_l\subseteq [i-1]\}$.
If $2r-q-1\le k$, then set
\[
\psi(\sigma)=\{B_1,\,\ldots,\,B_{q-1},\,B_q\cup B_r,\,\ldots,\,B_{r-1}\cup
B_{2r-q-1},\,B_{2r-q},\,\ldots,\,B_k\}.
\]
If $2r-q-1>k$, then set
\[
\psi(\sigma)=\{B_1,\,\ldots,\,B_{q-1},\,B_q
\cup B_r,\,\ldots,\,B_{q+k-r}
\cup B_k,\,B_{q+k-r+1},\,\ldots,\,B_{r-1}\}.
\]
\end{itemize}
\end{itemize}

It remains to show that the map
 $\psi$ is well-defined and it is indeed the inverse of
the map $\varphi$.

For any $\sigma\in \US_n\backslash\A_n$, we
notice that in Step 1 of the above construction of $\psi$,
the element $j$ always exists since $\sigma$ is not atomic.
Moreover, we observe that $j\geq 2$. In Step 2,
by the choice of $j$, we have
\begin{align*}
\sigma_{[j-1]}&=\{B_1,B_2,\ldots,B_{r-1}\},\\[5pt]
\sigma_{[j,\,n]}&=R(\sigma)=\{B_r,B_{r+1},\ldots,B_k\}.
\end{align*}
Since $\sigma$ is unsplitable,
we can always find the element $q$.
Otherwise, if every block $B_1,B_2,\ldots,B_k$ contains an element in
 $[i,n]$, by the assumption \eqref{j-1}, we have
$B_{p}\cap[i,n]\neq\emptyset$ for any $1\le p\le k$,
and
\[
\min(B_1\cap[i,n])
<\min(B_2\cap[i,n])
<\cdots<\min(B_k\cap[i,n]).
\]
This implies that
\[
\sigma=\sigma_{[i-1]}\circ(\sigma_{[i,n]}-i+1),
\]
a  contradiction to the fact that $\sigma$ is unsplitable.
This confirms the existence of the element $q$.
At this point, we still need to
show that $\psi(\sigma)\in\A_n\backslash\US_n$.
It is clear from
the above construction that $\psi(\sigma)$ is splitable.
For the case when $\sigma_{[j-1]}$ is
unsplitable,  it is easily seen
 that $\psi(\sigma)$ is atomic. When
$\sigma_{[j-1]}$ is splitable,
since $i\in B_1$ and $B_q\subseteq [i-1]$,  we find that $\psi(\sigma)$ is atomic.
Thus we have shown that $
\psi(\sigma)\in
\A_n\backslash\US_n$.
Consequently,
$\psi$ is well-defined.

It is not difficult to verify that $\psi$ is the inverse
map of $\varphi$. The details are omitted.
This completes the proof. \qed

The following example is an illustration of
the maps $\varphi$
and $\psi$. Let \[
\pi=\bigl\{\{1,3,5,6\},\,\{2,7,9\},\,\{4,8,10\}\bigr\}
\]
be the partition as given
in \eqref{eg_pi}.
In Step 1 of the map $\varphi$, we have $i=5$.
By \eqref{eg_raf}, we get
\begin{equation}
\label{vpi}
\varphi(\pi)=\{\{1,3,5,6\},\{2\},\{4\},\{7,9\},\{8,10\}\}.
\end{equation}

Conversely, assume that $\sigma$ is the
partition given in \eqref{vpi}.
Then $\psi(\sigma)$ is determined as follows.
First,  we have $R(\sigma)=\{\{7,9\},\,\{8,10\}\}$ and
$j=7$. Then,
\[ \sigma_{[j-1]}=\{\{1,3,5,6\},\,\{2\},\,\{4\}\}
=\{B_1,B_2,B_3\}
\]
 is
splitable, and $i=5$ is the smallest element
 in $\{1,3,5,6\}$ such
that
\[ \sigma_{[j-1]}=\sigma_{[i-1]}
\circ(\sigma_{[i,\,j-1]}-i+1).\]
Since $B_2$ is the first block of $\sigma_{[j-1]}$ that is
contained in $[i-1]$, we get
\[ \psi(\sigma)=\{\{1,3,5,6\},\,\{2,7,9\},\,\{4,8,10\}\}=\pi.\]

\noindent{\bf Acknowledgments.} This work was supported by the 973
Project, the PCSIRT Project of the Ministry of Science and
Technology, and the National Science Foundation of China.


\begin{thebibliography}{9}

\bibitem{BHRZ06}
N. Bergeron, C. Hohlweg, M.H. Rosas and M. Zabrocki, Grothendieck
bialgebras, partition lattices, and symmetric functions in
noncommutative variables, Electron. J. Combin. 13, 1 (2006),
\#R75.

\bibitem{BERGERON-R-R-Z}
N. Bergeron, C. Reutenauer, M.H. Rosas and M. Zabrocki, Invariants
and coinvariants of the symmetric group in noncommuting variables,
Canad. J. Math. 60, 2 (2008), 266--296.

\bibitem{BZ05X}
N. Bergeron and M. Zabrocki, The Hopf algebras of symmetric
functions and quasisymmetric functions in non-commutative
variables are free and cofree, arXiv:math.CO/0509265.

\bibitem{Can-Sagan}
M.B. Can and B.E. Sagan, Partitions, rooks, and symmetric
functions in noncommuting variables,
arXiv:math.CO/1008.2950.

\bibitem{Rosas-Sagan}
M.H. Rosas and B.E. Sagan,
Symmetric functions in noncommuting variables,
Trans. Amer. Math. Soc. 358, 1 (2006), 215--232.

\bibitem{Wolf}
M.C. Wolf, Symmetric functions of non-commutative elements,
Duke Math. J. 2 (1936), 626--637.

\end{thebibliography}
\end{document}